\documentclass[12pt]{amsart}

\usepackage{vmargin}
\setmargrb{1in}{1in}{1in}{1in}
\usepackage{amstext}
\usepackage{amssymb,mathrsfs}
\usepackage{latexsym}
\newcommand\RR{{\mathbb R}}

\newcommand\NN{{\mathbb N}}

\def\d={\,:=\,}

\newcommand{\semdir}
{\rtimes}

\font\frakten=eufm10
\newfam\frakfam
\textfont\frakfam=\frakten

\newtheorem{thm}{Theorem}
\newtheorem{lemma}[thm]{Lemma}
\newtheorem{cor}[thm]{Corollary}

\newtheorem{Defn}[thm]{Definition}
\newtheorem{Ex}[thm]{Example}
\newtheorem{Rem}[thm]{Remark}
\newtheorem{Exs}[thm]{Examples}
\newtheorem{Rems}[thm]{Remarks}
\newtheorem{Defrem}[thm]{Definition and Remark}
\newtheorem{Remnt}[thm]{}
\newenvironment{defn}
 {\begin{Defn} \begin{rm}} {\end{rm} \hfill $\Box$ \end{Defn}}

\newenvironment{rem}
 {\begin{Rem} \begin{rm}} {\end{rm} \hfill $\Box$ \end{Rem}}

\newenvironment{prf} {{\bf Proof.}}{\hfill $\Box$}

\begin{document}

\author{Hartmut F{\"u}hr}

\address{Lehrstuhl A f\"ur Mathematik \\
RWTH Aachen \\
52056 Aachen \\
Germany.}
%\address{Institute of Biomathematics and Biometry \\
%GSF Research Center for Environment and Health \\
%Ingolst{\"a}dter Stra{\ss}e~1 \\
%D-85764 Neuherberg \\
%Germany.}
%\email{fuehr@gsf.de}

\title[Generalized Calder\'on conditions]{Generalized Calder\'on conditions and regular orbit spaces}


\date{\today}

\begin{abstract}
The construction of generalized continuous wavelet transforms on locally
compact abelian groups $A$ from quasi-regular representations of a semidirect
product group $G = A \rtimes H$ acting on ${\rm L}^2(A)$
requires the existence of a square-integrable function whose Plancherel transform 
satisfies  Calder\'on-type resolution of the identity. The question then arises under what conditions 
such square-integrable functions exist.

The existing literature on this subject leaves a gap between sufficient and necessary criteria. 
In this paper, we give a characterization in terms of the natural action of the dilation
group $H$ on the character group of $A$.  
We first prove that a Calder\'on-type resolution of the identity gives rise to a decomposition of Plancherel
measure of $A$ into measures on the dual orbits, and then show that the latter property is equivalent to
regularity conditions on the orbit space of the dual action. 

Thus we obtain, for the first time, sharp necessary and sufficient criteria for the existence of a wavelet inversion formula.
As a byproduct and special case of our results we obtain that discrete series subrepresentations of the 
quasiregular representation correspond precisely to dual orbits with positive Plancherel measure and
associated compact stabilizers. Only sufficiency of the conditions was previously known. 
\end{abstract}

\maketitle

\section{Introduction}

The continuous wavelet transform of $f \in {\rm L}^2(\mathbb{R})$ is obtained by picking
a suitable $\psi \in {\rm L}^2(\mathbb{R})$ and letting 
\[
 V_\psi f (b,a) = \int_{\mathbb{R}} f(t) |a|^{-1/2} \overline{\psi \left( \frac{t-b}{a} \right)} dt~\mbox{ for } b  \in \mathbb{R}, a \in 
\mathbb{R} \setminus \{ 0 \}~.
\]
Among the many useful aspects of wavelets, probably the most fundamental one is wavelet inversion,
usually formulated as
\[
  f(t) = \int_{\mathbb{R}} \int_{\mathbb{R}} V_\psi f(b,a) |a|^{-1/2} \psi \left( \frac{t-b}{a} \right) db \frac{da}{a^2} ~,
\] to be read in the weak sense (rather than pointwise). This remarkable identity holds precisely if $\psi$ was chosen as 
{\bf admissible vector}, fulfilling the {\bf Calder\'on condition} 
\begin{equation}
\label{eqn:Cald_R1}
\int_{\mathbb{R}} \frac{|\widehat{\psi}(\xi)|^2}{|\xi|} d\xi = 1~.
\end{equation} 

The generalization of this construction, in particular to higher-dimensional euclidean space, 
has been studied early on, see e.g. \cite{Mu,Bo}. In the euclidean
setting, the role of the dilations $a \not=0$ is assumed by the elements of a matrix group $H$,
and various sources have studied which properties of $H$ ensure the existence of an inversion
formula, see e.g. \cite{BeTa,Fu,FuMa,LWWW}.
A further extension, replacing $\mathbb{R}^d$ by a general locally compact group $A$
and $H$ by a group of topological automorphisms, was considered in \cite{Ka_etal}.

The wavelet inversion formula is closely related to a suitable generalization of the Calder\'on
condition. As will be seen in the next section, this condition is quite easy to write down. However, it
is not at all trivial to decide whether there actually exist ${\rm L}^2$-functions satisfying it.  
Sufficient conditions for dilation groups acting on $\mathbb{R}^d$ were derived in \cite{FuMa,LWWW},
along with some necessary conditions. However, a complete characterization of these groups in terms
of necessary and sufficient conditions has been missing. The chief contribution of this paper is to
provide such a characterization in terms of the natural action on the dual group.

The paper is structured as follows: Section \ref{sect:wavelet_semdir} contains a more detailed exposition of the
group-theoretic construction of continuous wavelet transforms from the action of an automorphism
group on a locally compact abelian group.  We investigate wavelet inversion formulae valid 
for elements from a proper closed invariant subspace. For this purpose,  we introduce the dual action of 
the dilation group, and formulate the Calder\'on condition for admissible vectors.
A useful auxiliary notion for the discussion of admissible vectors is ``weak admissiblity''. 
We formulate a full characterization of dilation groups admitting weakly admissible vectors (Theorem \ref{thm:main_grp}),
which is the central result of this paper. 
The following two sections are devoted to a proof of this theorem. As it turns out, the core result
is of a predominantly measure-theoretic nature, and our treatment highlights these aspects. The
main result of these sections is Theorem \ref{thm:main}. 
In the final section we resume the discussion of admissible vectors. Theorems \ref{thm:main_grp}
and \ref{thm:main_ex_adm} provide a complete characterization of invariant subspaces allowing
a wavelet inversion formula. We also comment on irreducible subspaces with wavelet inversion formula,
which necessarily correspond to orbits of the dual action with positive measure
and compact fixed groups (Corollary \ref{cor:discrete_series}). 

\section{Wavelet transforms from semidirect products} \label{sect:wavelet_semdir}

Let us shortly sketch the group-theoretic framework for the construction of 
continuous wavelet transforms on locally compact abelian groups. The case
where the underlying group is $\mathbb{R}^n$ has been studied e.g. in \cite{BeTa,FuMa,
LWWW}, the generalization to arbitrary LCA groups was considered in \cite{Ka_etal}.

Let $A$ denote a second countable locally compact abelian group (with group structure
written additively), and let $H$ be
a group of topological automorphisms of $A$. $H$ denotes a group of automorphisms
of $A$, endowed with a secound countable locally compact
group topology making the natural action of $H$ on $A$ continuous. The
semidirect product group $A \rtimes H$ consists of elements $(a,h) \in A \times H$, with group law 
$(a,h)\cdot (b,g) = (a+ h(b),h g)$. When endowed with the product topology,
$G$ is a second countable locally comapct group as well. 

For any locally compact group $S$, integration  against (left) Haar measure is denoted as $\int_S g(s) ds$.
Haar measure of a Borel $B \subset S$ is denoted by 
$|B| = \int_S \mathbf{1}_B(s) ds$. Here, as below, we use the notation $\mathbf{1}_B$ for the indicator
function of $B$.  

The action of $H$ on $A$ induces a continuous homomorphism $\delta: H \to \mathbb{R}^+$
by $\delta(h) = \frac{|h(B)|}{|B|}$, where $ B \subset A$ is any Borel set of 
positive measure. The left Haar integral on $G$ is given by
\[
 \int_G f(x,h) d(x,h) = \int_H \int_A f(x,h) ~dx~\frac{dh}{\delta(h)}~,
\]
and the modular function of $G$ is $\Delta_G(a,h) = \frac{\Delta_H(h)}{\delta(h)}$. 

$G$ has a natural unitary representation acting on ${\rm L}^2(A)$ via
\[
  \pi(a,h) f(t) = \delta(h)^{-1/2} f(h^{-1}(t-a))~(t \in A)~.
\] Given a function $g \in {\rm L}^2(A)$,  the associated wavelet transform is an 
operator $V_\psi$ mapping $f \in {\rm L}^2(A)$ to its coefficient function $V_\psi f$, defined on $G$ as
\[
 V_\psi f (a,h) = \langle f,\pi(a,h) \psi \rangle ~. 
\]

\begin{defn} Let $\mathcal{H} \subset {\rm L}^2(A)$ be a closed $\pi$-invariant subspace. 
 $g \in \mathcal{H}$ is called {\bf weakly admissible (for $\mathcal{H})$} if $V_\psi : \mathcal{H} \to {\rm L}^2(G)$
 is a (well-defined) bounded injective
map. It is called {\em admissible (for $\mathcal{H})$} if $V_\psi : \mathcal{H} \to {\rm L}^2(G)$ 
is an isometric embedding. 
\end{defn}

\begin{rem}
 Admissibility is equivalent to a weak-sense inversion formula: $g$ is admissible for
$\mathcal{H}$
iff for all $f \in \mathcal{H}$,
\[
 f = \int_H \int_A V_\psi f(a,h)  \pi(a,h) \psi~ da \frac{dh}{\delta(h)}~,
\] holds in the weak sense (see e.g., \cite[Section 2.2]{Fu_LN}).
\end{rem}

\begin{rem}
Note that $V_\psi$ intertwines the action of $\pi$ with left translation. In particular, $g$
is weakly admissible iff $V_\psi$ is a bounded injective intertwining operator between the restriction
of $\pi$ to $\mathcal{H}$ and the left regular representation acting on ${\rm L}^2(G)$. 
In fact, the existence of weakly admissible is equivalent to unitary containment in the regular representation \cite[2.21]{Fu_LN}.
\end{rem}

We denote the dual group of $A$ by $\widehat{A}$. It is a second countable locally compact abelian group
as well. (For this and the following facts concerning locally compact abelian groups, see \cite{Fo}.) 
The Fourier transform of $f \in {\rm L}^1(A)$ is defined as
\[
 \widehat{f}(\xi) = \int_A f(x) \overline{\xi(x)} dx~.
\] We normalize Haar measure on $\widehat{A}$ such that for all $f\in {\rm L}^1(G) \cap {\rm L}^2(G)$,
$\| f \|_2 = \| \widehat{f} \|_2$. The Plancherel theorem implies that the Fourier transform
extends to a unitary operator ${\rm L}^2(A) \to {\rm L}^2(\widehat{A})$. For this reason, Haar
measure on $\widehat{A}$ is also called Plancherel measure of $A$.

The action of $H$ on $A$ gives rise to the dual action on $\widehat{A}$, which is
a right action defined by $(\xi.h)(x) = \xi(h(x))$. The behaviour of Haar measure on $\widehat{A}$
is similar to that of Haar measure on $A$, i.e., $|B.h| = \delta(h) |B|$ for all $B \subset \widehat{A}$
Borel. 

For the study of (weakly) admissible
vectors for invariant subspaces, the dual action is an indispensable tool. To begin with,
invariant subspaces are in one-to-one correspondence to $H$-invariant Borel subsets of 
$\widehat{A}$, by the following result.
\begin{lemma} \label{lem:char_inv_subsp}
 Let $X \subset \widehat{A}$ be an $H$-invariant Borel subset. Let 
\[
 \mathcal{H}_X = \{ f \in {\rm L}^2(A) : \widehat{f} \cdot \mathbf{1}_X = \widehat{f} \}~.
\] Then $\mathcal{H}_X \subset {\rm L}^2(A)$ is a $\pi$-invariant closed subspace.
We write $\pi_X$ for the restriction of $\pi$ to $\mathcal{H}_X$.

Conversely, if $\mathcal{H}  \subset  {\rm L}^2(A)$ is a $\pi$-invariant and closed subspace,
then $\mathcal{H} = \mathcal{H}_X$ for a suitably chosen $H$-invariant Borel set $X$.
\end{lemma}
\begin{prf}
 First note that if $\mathcal{H}$ is invariant under shifts, i.e. all operators of the type 
$\pi(a,e_H)$, then necessarily $\mathcal{H} = \mathcal{H}_X$ for some Borel set $X$.
This follows from the characterization of the commuting algebra by the Fourier transform,
e.g. \cite[4.44]{Fo}. If, in addition, $\mathcal{H}$ is also invariant under $\pi (0, h)$
for all $h$, it necessarily follows that, possibly after removing a set of measure zero, $X$ is 
in addition $H$-invariant.  The proof given in \cite{Fu} for this fact in the case $A = \mathbb{R}^n$ carries over verbatim.  
\end{prf}

We next turn to the derivation of admissibility criteria. Direct calculation employing
the Plancherel Theorem for $A$ allows to derive the crucial equality
\[
 \| V_\psi f \|_2^2 = \int_{\widehat{A}} | \widehat{f} (\xi) |^2 \int_{H} |\widehat{\psi}(\xi.h)|^2 ~dh~d\xi~.
\] See e.g. \cite{FuMa,LWWW} for the proof in the case $A= \mathbb{R}^d$, which immediately carries over
to the general setting.
From this, one easily derives the following criteria for strong and weak admissibility, generalizing the Calder\'on
condition for wavelets over the reals:  
\begin{lemma} \label{lem:calderon}
Let $\mathcal{H} \subset {\rm L}^2(A)$ be closed and $\pi$-invariant. Hence $\mathcal{H}
= \mathcal{H}_X$ for a suitable $H$-invariant Borel set $X \subset \widehat{A}$. Then
 $\psi \in \mathcal{H}$ is weakly admissible iff the function 
\begin{equation} \label{eqn:cald_fn}
 \xi \mapsto \int_H |\widehat{\psi}(\xi.h)|^2 ~dh
\end{equation} is a.e. bounded and nonzero on $X$. Moreover,
$\psi$ is admissible iff this function equals one a.e. 
\end{lemma}

Furthermore, it is easily verified that for
\[
 \int_H |\widehat{\psi}(\xi.h)|^2 ~dh
\] to be finite, the stabilizer of $\xi$, defined by $H_\xi =  \{ h \in H : \xi.h = \xi \}$
must be compact (see Lemma \ref{lem:cmp_stab} below).  
Hence Lemma \ref{lem:calderon} implies that almost all stabilizers must be compact for $H$ to 
be weakly admissible. 
However, it has been noticed early on that this necessary condition is not sufficient:
The relevant counterexample is provided by letting $A = \mathbb{R}^2$ and $H = {\rm SL}(2,\mathbb{Z})$.
It turns out that almost all stabilizers are trivial, but $H$ is not weakly admissible (see \cite{Fu} for a related
example).

Additional sufficient criteria were provided in \cite{FuMa,Fu_LN,LWWW} for the case of $A = \mathbb{R}^n$
and a matrix group $H$, but the results in these
papers do not yield a full characterization.
The authors of \cite{LWWW} studied the condition that for almost every $\xi$ there exists $\epsilon>0$
such that the $\epsilon$-stabilizer $H_{\epsilon,\xi} = \{ h \in H : |\xi - \xi.h|< \epsilon \}$ is a compact subset
of $H$. Here $| \cdot |$ denotes the Euclidean distance. It is shown in \cite{LWWW} that this
condition ensures weak admissibility. Necessity of this condition was conjectured, but not shown in \cite{LWWW}.
By contrast, \cite{FuMa,Fu_LN} studied regularity conditions on the orbit spaces, somewhat similar to the 
properties that will be considered in the next section. However, no necessary condition was derived. 

The following theorem is the chief result of this paper. It characterizes the groups $H$ allowing a weakly
admissible vector in terms of regularity properties of the orbit space. 
\begin{thm} \label{thm:main_grp} Let $\mathcal{H} =\mathcal{H}_X$, for $X \subset \widehat{A}$
Borel and $H$-invariant. 
 $\mathcal{H}$ has a weakly admissible vector iff there exists a conull $H$-invariant Borel subset $B \subset X$ such that
\begin{enumerate}
 \item For all $\xi \in B$, the stabilizer $H_\xi$ is compact.
 \item There exists a Borel set $C \subset B$ such that for all $\xi \in B$, the set $C \cap \xi.H$ is
a singleton. 
\end{enumerate}
\end{thm}
This result is a direct consequence of the purely measure-theoretic Theorem 
\ref{thm:main} below. We have chosen to remove (almost) all references to wavelets and
harmonic analysis from the following two sections, because we believe that the central
problem is measure-theoretic in nature, and of a certain independent interest.

\section{Measure-theoretic setup and main result}  \label{sect:meastheory}

Let us begin by fixing terminology. A useful survey of the relevant definitions and results 
concerning Borel spaces can be found in \cite{AuMo}.

A Borel space is a set $X$ endowed with a $\sigma$-algebra on $X$.
The elements of the $\sigma$-algebra are called Borel sets. A measure defined on the $\sigma$-algebra is called Borel measure. 
A map between Borel spaces is called Borel if the preimage of Borel sets are Borel again.
A Borel isomorphism is a bijection between Borel spaces that is Borel in both directions. 
In the following, the Borel structures of locally compact groups and metric spaces 
are understood to be generated by the respective topologies. 
A Borel space is called standard if it is Borel isomorphic to a Borel subset of a separable 
complete metric space. We note that second countable locally compact groups are completely
metrizable, and therefore standard. This applies to $H$, but also to $A$ and
$\widehat{A}$.  Also, Borel subsets of standard spaces are clearly standard. 
Throughout the following two sections, $X$ denotes a
standard Borel measure space, on which a locally compact second countable group $H$ acts jointly 
measurably from the right.

We assume to be given a fixed $\sigma$-finite Borel measure $\lambda$ on 
$X$, which is quasi-invariant under $G$. This means that for all $h \in H$, the measure
$\lambda_h: A \mapsto \lambda(A.h)$ is equivalent to $\lambda$. By the Radon-Nikodym theorem, this
assumption implies the existence of a function $\rho: X \times H \to \mathbb{R}$ such that 
\[
 \frac{d\lambda_h}{d\lambda}(\xi) = e^{\rho(\xi,h)}
\] holds, for all $(\xi,h) \in X \times H$. $\rho$ is called the cocycle of the measure; it can be assumed
measurable on $X \times H$, and such that it fulfills the following cocycle conditions, for all $g,h \in H$ and $\xi \in  X$: 
\begin{eqnarray} \label{eqn:cc_cond}
 \rho(\xi,gh) & = & \rho(\xi.g,h) + \rho(\xi,g)~, \\
 \rho(\xi,h) & = & 0 \mbox{ if } \xi.h =\xi~, \label{eqn:cc_fix}
\end{eqnarray} see e.g. \cite{GrSch} or \cite[Appendix B]{Zi}. 
We note that the definition of the cocycle entails the following two formulae for integration:
\begin{eqnarray}
  \lambda(B.h) & = & \int_B e^{\rho(\xi,h)} d\lambda(\xi) \label{eqn:trans_form_meas} \\
 \int_X f(\xi.h^{-1}) d\lambda(\xi) & = & \int_X f(\xi) e^{\rho(\xi,h)} d\lambda(\xi) \label{eqn:trans_form_fun}~,
\end{eqnarray} where the second equation holds for all positive Borel functions $f$, in the extended sense
that one side is infinite iff the other is. 

%  We will call the measure $\lambda$ {\bf tame}
% if it has a cocycle $\rho$ satisfying the condition 
% \begin{equation}
%  \label{eqn:tame_cc}
%  \forall \xi \in  X \forall C \subset H \mbox{ compact }: \rho(\xi,C) \subset \mathbb{R} \mbox{ is bounded.}
% \end{equation} 

We denote by $X/H$ the space of all $H$-orbits in $X$. Let $q: X \to X/H$
denote the quotient map. $X/H$ is endowed
with the quotient Borel structure: A subset $B \subset X/H$ is declared
Borel if $q^{-1}(B) = \bigcup \{ W : W \in B \} \subset X$ is Borel. 
% 
% As a consequence of Tonelli's theorem, for any Borel function $\varphi : X
% \to \mathbb{R}_0^+$, letting
% \[
% g(\xi) = \int_H f(\xi.h) dh~
% \] defines a Borel map $g: X \to \mathbb{R}_0^+ \cup \{ \infty \}$.
% Moreover, since the action is on the right, and Haar measure is
% left-invariant, $g$ is $H$-invariant, �and can thus be regarded as a
% Borel function on $X/H$.

\begin{defn}
The action of $H$ on $X$ is called {\bf weakly admissible} if there
exists a Borel function $\varphi : X \to \mathbb{R}^+$ satisfying
 \begin{equation} \label{eqn:weak_adm_fun}
0 <  \int_H \varphi(\xi.h) dh < \infty~, \mbox{ for $\lambda$-a.e. } \xi \in  X~.
\end{equation}
\end{defn}

The definition is a clear analogy to the Calder\'on condition. If
$X \subset \widehat{A}$ is an invariant Borel subset, we will see shortly 
that the existence of weakly admissible vectors for
a representation
is equivalent to weak admissibility of the dual action. The following
lemma spells out the technical details. 
\begin{lemma}
If the action of $H$ is weakly admissible, there exists a function $\varphi$ 
such that 
\begin{equation}
 \label{eqn:wCalderon_L1}
\varphi \ge 0~,~0 < \int_H \varphi(\xi.h)dh \le 1~(\lambda\mbox{-a.e})~,~\varphi \in {\rm L^1}(X,\lambda)~.
\end{equation}
\end{lemma}
\begin{prf}
Write $X$ as a disjoint union of sets of finite measure, $X = \bigcup_{n \in \mathbb{N}} X_n$.
Assume that $\varphi_1$ fulfills (\ref{eqn:weak_adm_fun}), and define 
\[ \varphi_2(\xi) = \frac{\min(1,\varphi_1(\xi))}{2^n (1+\lambda(X_n))} ~,~\xi \in  X_n ~.\]
Then $\varphi_2$ is integrable, and also fulfills (\ref{eqn:weak_adm_fun}). The same is then true for 
\[
 \varphi(\xi) = \frac{\varphi_2(\xi)}{1+\int_H\varphi_2(\xi.h) dh}~,
\] which is the desired function.  
\end{prf}

\begin{cor} \label{cor:equiv_wadm_exwadmvec}
Let $\mathcal{H} = \mathcal{H}_X\subset {\rm L}^2(A)$, for a suitable $H$-invariant $X \subset\widehat{A}$.
There exists a weakly admissible vector for $\mathcal{H}$ iff the dual action of $H$ on $\widehat{X}$ is weakly admissible.
\end{cor}
\begin{prf}
The ``only-if''-part is clear. For the other direction, 
we let $\widehat{\psi}(\xi) = \varphi(\xi)^{1/2}$, where $\varphi$ fulfills (\ref{eqn:wCalderon_L1}). 
\end{prf}

As already indicated in the title, the structure of the orbit space $X/H$ is of
central interest. Such spaces can be quite pathological.  By contrast, the situation for each
individual orbit is quite simple, as the following lemma shows.
\begin{lemma}\label{lem:Borel_ind_orbits}
 For all $\xi \in  X$, the orbit $\xi. H \subset X$ is Borel. 
Furthermore, the stabilizer $H_\xi = \{ h \in H : \xi.h = \xi \}$
is a closed subgroup of $H$, and the quotient map $H \ni h \mapsto \xi.h$
induces a Borel isomorphism $H_\xi \setminus H \to \xi.H$. 
\end{lemma}
\begin{prf}
Confer \cite[Chapter I, Proposition 3.7]{AuMo}. 
\end{prf}

Now the first necessary condition for weakly admissible actions is easily
proved.  
\begin{lemma} \label{lem:cmp_stab}
 The set 
\[
 X_c := \{ \xi \in X : H_\xi \mbox{ is compact } \} \subset X
\] is Borel and $H$-invariant. If $H$ is weakly admissible, then $X_c \subset X$ is
conull. 
\end{lemma}
\begin{prf}
 The stabilizer map $x \mapsto H_\xi$ is Borel,
if one endows the set of closed subgroups with the Fell topology \cite[Chapter II, Proposition 2.3]{AuMo}.
Moreover, the set of compact subgroups is Borel (see \cite[Proposition 5.5]{Fu_LN}), hence $X_c$ is Borel. 
$H$-invariance is immediate
from the observation that all stabilizers associated to a given orbit are conjugate. 

If $\varphi$ is a positive Borel function on $X$ and $\xi \in X$
is such that 
\[
 0 < \int_H \varphi(\xi.h) dh < \infty
\] then the fact that the function $h \mapsto \varphi(\xi.h)$ is integrable (with nonzero integral)
and leftinvariant under the closed subgroup $H_\xi$ at the same time forces $H_\xi$ to be compact:
E.g., pick $\epsilon>0$ such that $|\{ h : \varphi(\xi.h) > \epsilon \}|>0$. This set has finite
Haar measure, and is left $H_\xi$-invariant, thus \cite[Lemma 11]{Fu} yields compactness of $H_\xi$.

Thus, if the action of $H$ is weakly admissible, $|X \setminus X_c|=0$. 
\end{prf}

We will characterize weak admissibility in terms of
measure-theoretic properties of $X/H$, which are closely
related to standardness.  A Borel space is called {\bf countably generated} if
the $\sigma$-algebra is generated by a countable subset. 
It is called {\bf separated} if single points are Borel.
A Borel space is called {\bf countably separated} if there is a sequence of
Borel sets separating the points. All these properties are
inherited by products and Borel subspaces.
A Borel space is called {\bf analytic} if it is (Borel-isomorphic to) the
Borel image of a standard space in a countably generated space.

We say that $X/H$ admits a $\lambda$-transversal if there exists an $H$-invariant
$\lambda$-conull Borel set $Y \subset X$ and a Borel set $C\subset Y$ meeting each orbit 
in $Y$ in precisely one point.

A {\bf pseudo-image} of $\lambda$ is a measure $\overline{\lambda}$
on $X/H$ obtained as image measure of an equivalent finite measure under the quotient map $q$;
clearly all pseudo-images are equivalent.  We call
$\overline{\lambda}$ {\bf standard}, if there exists $Y \subset X$ Borel,
$H$-invariant, conull, such that $Y/H$ is standard.

Finally, we need the notion of a {\bf measure decomposition}: 
A {\bf measurable family of measures} is a family $(\beta_{\mathcal{O}})_{\mathcal{O}
\subset X}$ indexed by the orbits in $X$, such that for all Borel sets $B \subset X$,
the map $\mathcal{O} \mapsto \beta_{\mathcal{O}}(B)$ is Borel on $X/H$. 

A  {\bf measure decomposition of $\mathbf{\lambda}$} consists of a pair 
$(\overline{\lambda},(\beta_{\mathcal{O}})_{\mathcal{O}
\subset X})$ , where  $\overline{\lambda}$
is a pseudo-image of $\lambda$ on $X/H$, or a $\sigma$-finite measure equivalent
to such a pseudo-image, and a measurable family $(\beta_{\mathcal{O}})_{\mathcal{O}
\subset X}$ such that for all $B \subset X$ Borel,
\[
 \lambda(B) = \int_{X/H} \beta_{\mathcal{O}}(B) d\overline{\lambda}(\mathcal{O})~.
\]
Note that this entails, for all positive Borel functions $f$ on $X$, that
\[
 \int_{X} f(\xi) d\lambda(\xi) = \int_{X/H} \int_{\mathcal{O}} f(\xi) d\beta_{\mathcal{O}}(\xi) d\overline{\lambda}(\mathcal{O}). 
\]
We say that {\bf $\mathbf{\lambda}$ decomposes over the orbits} if there exists a measure decomposition
with the additional requirement that, for $\overline{\lambda}$-almost every $\mathcal{O} \in  X/H$, 
the measure $\beta_{\mathcal{O}}$ is supported in $\mathcal{O}$, meaning 
$\beta_{\mathcal{O}}(X\setminus \mathcal{O}) = 0$. 

\begin{thm} \label{thm:main}
Let $X$ be a standard Borel-space, and $H$ a second-countable group acting 
measurably on $X$.  Assume that $\lambda$ is a quasi-invariant $\sigma$-finite measure on $X$.
 Consider the following statements:
\begin{enumerate}
\item[(a)] The action of $H$ is weakly admissible.
\item[(b)] $\lambda$ decomposes over the orbits. 
\item[(c)] $\overline{\lambda}$ is standard.
% \item[(d)] $X/H$ admits a measurable cross-section.
\item[(d)] $X/H$ admits a $\lambda$-transversal.
\end{enumerate}
Then $(a) \Rightarrow (b) \Leftrightarrow (c) \Leftrightarrow (d)$, and $\lambda(X \setminus X_c) = 0$. 

Conversely, if $\lambda(X \setminus X_c) = 0$, then $(d) \Rightarrow (a)$. 
\end{thm}

The equivalence of $(b)$ through $(d)$ is possibly folklore,
although we have not been able to locate a handy reference for the "almost-everywhere" version that we
consider here. Also, the proof of $(b) \Rightarrow (c)$ turned out to be rather more technical than
initially expected. We include detailed arguments for the sake of reference.

Note that Corollary \ref{cor:equiv_wadm_exwadmvec} and Theorem \ref{thm:main}, applied to the dual action of $H$
on the invariant set $X$, indeed imply Theorem \ref{thm:main_grp}. 

\section{Proof of Theorem \ref{thm:main}}

If the action of $H$ is weakly admissible, then $\lambda(X\setminus X_c)=0$ by Lemma \ref{lem:cmp_stab}.
W.l.o.g., we will therefore assume in the following that $X=X_c$.

\subsection{Proof of $(a) \Rightarrow (b)$}

Assume that the action of $H$ is weakly admissible. 
The proof strategy will be to define a measure $\mu$ that decomposes over the orbits, 
and to show that $\mu$ is $\sigma$-finite and equivalent to $\lambda$. 

\begin{lemma} \label{lem:ovrlmbd_phi}
 Let $\overline{\lambda}$ be a pseudo-image of $\lambda$ on $X/H$. Let  $\varphi: X \to \RR_0^+$ be a
 Borel function satisfying (\ref{eqn:wCalderon_L1}). For Borel sets
 $U \subset X/H$ let
 \[
  \overline{\lambda_\varphi} (U) = \int_{q^{-1}(U)}
  \varphi(\xi) ~d\lambda(\xi)~~.
 \] Then $\overline{\lambda_\varphi}$ is a finite measure on
 $X$ that is equivalent to $\overline{\lambda}$, satisfying for all Borel functions
 $g: X/H \to \mathbb{R}_0^+$
\begin{equation} \label{eqn:int_lphi}
 \int_{X/H} g(\mathcal{O}) d\overline{\lambda}_\varphi = \int_X g(\xi.H) \varphi(\xi) d\lambda(\xi)~.
\end{equation}
\end{lemma}

\begin{prf}
 We need to show for an $H$-invariant Borel set $V \subset
 X$ that
 \[ \lambda(V) = 0 \Leftrightarrow \int_V \varphi
 (\xi) d\lambda(\xi) = 0~~. \]
 Direction ''$\Rightarrow$'' is clear, since the right-hand side is an integral
 over a $\lambda$-null set.

 For the other direction, we employ the quasi-invariance of
 $\lambda$ and invariance of $V$ to note that for all $h \in H$,
 \[
 \int_V \varphi(\xi.h) d\lambda(\xi) = 0 ~.
 \] Integrating over $H$ and applying Tonelli's theorem, we obtain
 \[
 \int_V \int_H \varphi(\xi.h) ~dh~ d\lambda(\xi) = 0~.
 \] By assumption, the inner integral vanishes $\lambda$-almost
 nowhere, hence $\lambda(V) = 0$ follows.

By definition, equation (\ref{eqn:int_lphi}) holds for indicator functions, and extends to 
nonnegative Borel maps by standard arguments. 
\end{prf}

%From now on, we assume $\overline{\lambda}_\varphi = \overline{\lambda}$.
\begin{lemma}
\label{lem:defnmu} 
\begin{enumerate}
 \item[(a)] Let $\mathcal{O} = \xi.H$, and assume that $H_\xi$ is compact. Then
\[
 \mu_{\mathcal{O}}(B) = |\{ h \in H : \xi . h \in B \}|
\] defines a $\sigma$-finite measure supported on $\mathcal{O}$.
$\mu_{\mathcal{O}}$ is independent of the choice of $\xi \in \mathcal{O}$. 
\item[(b)] Let $\varphi$ be a Borel function satisfying (\ref{eqn:wCalderon_L1}).  
For Borel sets $B \subset X$, define
\begin{eqnarray} \label{eqn:mu}
 \mu(B) & = & \int_{X/H} \mu_{\mathcal{O}}(B)
 ~d\overline{\lambda}_\varphi(\mathcal{O}) \\ \nonumber
 & = & \int_{X} \varphi(\xi) \int_H
 \mathbf{1}_B(\xi.h) ~dh d\xi ~~.
\end{eqnarray}
$\mu$ is a well-defined Borel measure.
\end{enumerate}
\end{lemma}

\begin{prf}
Since $\mathcal{O}$ is Borel, $\mu_{\mathcal{O}}$ is a well-defined Borel-measure. Furthermore, since $H_\xi$ is compact,
$\mu_{\mathcal{O}}$ is finite on sets of the form $\xi.C$, with $C$ compact. In particular, 
since $H$ is $\sigma$-compact, $\mu_{\mathcal{O}}$
is $\sigma$-finite. $\mu_{\mathcal{O}}$ is independent of the choice of $\xi$, since
the action is on the right, and Haar-measure on $H$ is leftinvariant. 
The well-definedness of $\mu$ follows
 from Fubini's theorem and the measurability of $(h,\xi)
 \mapsto \mathbf{1}_A(\xi.h)$. 
 The second equation of (\ref{eqn:mu}) is obtained directly from (\ref{eqn:int_lphi}).
\end{prf}

The following result will allow to establish equivalence $\mu$ and $\lambda$.
\begin{lemma} \label{lem:phi_dec_lambda}
Let $\varphi$ be a positive Borel function fulfilling
(\ref{eqn:wCalderon_L1}), and let \[ \Phi(\xi) = \int_{H} \varphi(\xi.h)
~dh ~~.\] Then, for all Borel functions $f : X \to \RR_0^+$,
\begin{equation} \label{eqn:phi_dec_lambda}
\int_{X} f(\xi) d\lambda(\xi) = \int_X \frac{\varphi(\xi)}{\Phi(\xi)} \int_H f(\xi.h) e^{\rho(\xi,h)} \Delta_H(h)^{-1} ~dh ~d\lambda(\xi) ~~.
\end{equation}
\end{lemma}

\begin{prf}
 The proof is a straightforward computation, using Tonelli's theorem:
 \begin{eqnarray*}
  \int_X f(\xi) ~d\lambda(\xi) & = & \int_X
  \frac{f(\xi)}{\Phi(\xi)}
  \int_H \varphi(\xi.h) ~dh~d\lambda(\xi) \\
  & = & \int_H \int_X \frac{f(\xi)}{\Phi(\xi)} \varphi(\xi.h) ~d\lambda(\xi)~dh \\
 & = & \int_H \int_X \frac{f(\xi)}{\Phi(\xi)} \varphi(\xi.h) ~d\lambda(\xi) ~dh \\
& = & \int_H \int_X \frac{f(\xi.h^{-1})}{\Phi(\xi.h^{-1})} \varphi(\xi) e^{\rho(\xi,h^{-1})} ~d\lambda(\xi) ~dh \\
& = & \int_X \frac{\varphi(\xi)}{\Phi(\xi)} \int_H f(\xi.h^{-1}) e^{\rho(\xi,h^{-1})} ~dh ~d\lambda(\xi) \\
 & = &  \int_X \frac{\varphi(\xi)}{\Phi(\xi)} \int_H f(\xi.h) e^{\rho(\xi,h)} \Delta_H(h)^{-1} ~dh ~,
 \end{eqnarray*} where the penultimate equality used $H$-invariance of $\Phi$. 
\end{prf}
% 
% It is tempting to read (\ref{eqn:phi_dec_lambda}) as a measure
% decomposition. Note however that the map
% \[
% \xi \mapsto \frac{\varphi(\xi)}{\Phi(\xi)} \int_H f(\xi.h^{-1}) e^{\rho(\xi,h^{-1})} ~dh
% \] usually varies along a single orbit, hence cannot be understood as the result of
% integrating $f$ over $\mathcal{O}$ against a fixed measure depending
% only on $\mathcal{O}$.

The next lemma establishes $\sigma$-finiteness of $\mu$:
\begin{lemma} \label{lem:mu_sf}
Let $H$ be weakly admissible.
\begin{enumerate}
\item[(a)] There exists $\varphi: X \to \RR_0^+$ satisfying (\ref{eqn:wCalderon_L1}), 
and in addition, the map $h \mapsto \varphi(\xi.h)$ is
continuous, for all $\xi \in  \Omega_\varphi$.
 \item[(b)] Let $\varphi$ satisfy (\ref{eqn:cald_fn}). For $k \in\NN$, define
 \[
  A_k = \{ \xi \in  X : \varphi(\xi) > 1/k \}~.
 \]
Then, for all $k \in \NN$ and $g \in H$:
\begin{equation}
 \mu(A_k.g) = \Delta_H(g) \mu(A_k) \le \Delta_H(g) k
 \overline{\lambda_\varphi}(X/H) < \infty~.
\end{equation}
\item[(c)] With $\varphi, \Omega_\varphi$ as in part (a), and
$A_k$ as in part (b): If $(h_n)_{n \in \NN} \subset H$ is dense,
then $\Omega_\varphi \subset \bigcup_{n,k \in \NN} A_k.h_n$. \\
In particular, $\mu$ is $\sigma$-finite. 
\end{enumerate}
\end{lemma}
\begin{prf}
For the proof of $(a)$ pick $\varphi_0$ 
satisfying (\ref{eqn:wCalderon_L1}). Pick a continuous, compactly
supported $\nu: H \to \RR_0^+$ satisfying
\[
 \int_H \nu(g) \Delta_H(g) ~dg = 1~~.
\] Letting
\[
\Omega_\varphi = \{ \xi \in  X : \int_H \varphi_0(\xi.h) ~dh < \infty
\}~~
\]
defines an $H$-invariant conull Borel subset. For $\xi \in  \Omega_\varphi$, we define
\[
 \varphi(\xi) = \int_H \varphi_0(\xi.g) \nu(g) ~dg~,
\]
and obtain
 \begin{eqnarray*}
\varphi(\xi.h) & = & \int_H \varphi_0(\xi.hg) \nu(g) dg ~.
\end{eqnarray*} The assumption $\xi \in  \Omega_\varphi$ amounts
to saying that the map $g \mapsto \varphi_0(\xi.g)$ is in ${\rm
L}^1(H)$. Now strong continuity of the left action of $H$ on ${\rm
L}^1(H)$ and boundedness of $\psi$ imply that $h \mapsto \varphi (\xi.h)$ is
continuous. 

Integrability of $\varphi$ is a straightforward consequence of $\varphi_0 \in 
{\rm L}^1(X,\lambda)$, $\nu \in {\rm L}^1(H)$ and Fubini's theorem.  
Finally,
\begin{eqnarray*}
 \int_H \varphi(\xi.h) dh & = & \int_H \int_H
 \varphi_0(\xi.hg) \nu(g) ~dg~ dh \\
  & = & \int_H \nu(g) \int_H \varphi_0(\xi.hg) ~dh~dg \\
  & = & \int_H \nu(g) \Delta_H(g) \int_H \varphi_0(\xi.h) ~dh~dg
  \\ & = & \int_H \varphi_0(\xi.h) ~dh~~,
\end{eqnarray*}
where the last equation was due to our choice of $\nu$. Hence
(\ref{eqn:wCalderon_L1}) for $\varphi_0$ implies the same for
$\varphi$, and $(a)$ is shown.

The first equation of part $(b)$ follows from
\begin{eqnarray*}
 \int_H \mathbf{1}_{A_k.g}(\xi.h) ~dh & = & \int_H \mathbf{1}_{A_k}
 (\xi.hg^{-1})~dh \\
  & = & \Delta_H(g) \int_H \mathbf{1}_{A_k}(\xi.h)~dh~~,
\end{eqnarray*}
and integration over $X/H$. For the inequality, observe that by
definition of $A_k$, we have $\mathbf{1}_{A_k}(\xi) < k \varphi(\xi)$,
and thus by choice of $\varphi$
\[
 \mu_{\xi.H}(A_k) = \int_H \mathbf{1}_{A_k}(\xi.h) dh \le k \int_H \varphi(\xi.h)~dh \le k  ~.
\]
But then
\[ \mu(A_k) = \int_{X/H} \mu_{\mathcal{O}}(A_k) d\overline{\lambda}_{\varphi} (\mathcal{O})
 \le  k \overline{\lambda_\varphi
}(X/H) ~~,\]

For part $(c)$ let $\xi \in  \Omega_\varphi$, hence $0 < \int_H
\varphi(\xi.h) d\lambda(\xi) \le 1$. Hence the integrand cannot be
identically zero, and there exists $k \in \NN$ such that
\[ B = \{ g \in H : \varphi(\xi.g^{-1}) > 1/k \} \] is
nonempty. By choice of $\varphi$, $B$ is open, hence there exists $n
\in \NN$ such that $h_n \in B$, implying $\varphi(\xi.h_n^{-1})> 1/k$.
But this means $\xi \in  A_k.h_n$, as desired.
\end{prf}

Now the implication $(a) \Rightarrow (b)$ is easily proved. 
We pick $\varphi$ according to Lemma \ref{lem:mu_sf} (a), and consider the measure $\mu$
defined in Lemma \ref{lem:defnmu}, using  $\overline{\lambda} = \overline{\lambda}_\varphi$. 
Then $\mu$ is equivalent to $\lambda$: 
On the one hand, Lemma \ref{lem:phi_dec_lambda}
provides for an arbitrary Borel set $A \subset X$ 
\begin{equation} \label{eqn:lambdaA}
 \lambda(A) =\int_X \frac{\varphi(\xi)}{\Phi(\xi)} \int_H \mathbf{1}_A(\xi.h)  e^{\rho(\xi,h)} \Delta_H(h)^{-1} ~dh ~d\lambda(\xi)~~,
\end{equation} whereas by Lemma \ref{lem:ovrlmbd_phi},
\begin{equation} \label{eqn:muA}
\mu(A) = \int_{X} \varphi(\xi) \int_H
\mathbf{1}_A(\xi.h) ~dh ~d\lambda(\xi)~~.
\end{equation}
Hence, by (\ref{eqn:lambdaA}), $\lambda(A) = 0$ iff $
\varphi(\xi) \int_H \mathbf{1}_A(\xi.h)  e^{\rho(\xi,h)} \Delta_H(h)^{-1} ~dh = 0$ for $\lambda$-a.e. $\xi$. 
  Both $\Delta_H$ and the exponential function are strictly positive, hence this is the case precisely
 when
 $\varphi(\xi) \int_H \mathbf{1}_A(\xi.h)
 ~dh = 0$ for $\lambda$-a.e. $\xi$. But by (\ref{eqn:muA}), the
 latter case is equivalent to $\mu(A) = 0$. Hence $\lambda$ and
 $\mu$ are equivalent.

Recall that by definition, $d\mu(\xi) = d\mu_{\mathcal{O}}(\xi) d\overline{\lambda}(\mathcal{O})$.
By Lemma \ref{lem:mu_sf} (c), $\mu$ is $\sigma$-finite. Hence the Radon-Nikodym Theorem applies, and yields
\[ d\lambda(\xi) = \frac{d\lambda}{d\mu}(\xi) d\mu(\xi) =  \frac{d\lambda}{d\mu}(\xi) d\mu_{\mathcal{O}}(\xi) d\overline{\lambda}(\mathcal{O})~,
\] which shows that letting $d\beta_{\mathcal{O}} (\xi) = \frac{d\lambda}{d\mu}(\xi) d\mu_{\mathcal{O}}(\xi)$ yields the desired measure decomposition.

\subsection{Proof of $(b) \Rightarrow (c)$}

% Assume that we are given a measure decomposition
% $(\overline{\lambda}, (\beta_{\mathcal{O}})_{\mathcal{O} \subset
% X})$. $\sigma$-finiteness of $\mu$ necessarily forces
% $\sigma$-finiteness of the $\beta_{\mathcal{O}}$, at least almost
% everywhere. In addition, the $\beta_{\mathcal{O}}$ have the same 
% cocycle $\rho$ as $\lambda$.

For this step, we first replace $\lambda$ by an equivalent probability measure $\alpha$. 
Then $\alpha$ decomposes over the orbits as well, by the same argument as in the proof
of $(a) \Rightarrow (b)$.  
In the decomposition of $\alpha$, almost every $\beta_{\mathcal{O}}$ is finite, and can thus be 
normalized to be a probability measure. Then the measure on the quotient space
effecting the decomposition of $\alpha$ into the normalized measures turns out to be a probability measure as well.

In short, $\lambda$ can be assumed to be a probability measure, and all measures involved
in the decomposition as well. Furthermore, we may assume that $\overline{\lambda}$ is
the image measure of $\lambda$ under $q$. The following argument
relates the decomposition to the ergodic decomposition constructed in \cite{GrSch},
and then uses properties of the latter.
For this purpose, let $\rho$ denote the cocycle of $\lambda$.
Let $M_\rho(X)$ denote the set of Borel probability
measures on $X$ with cocycle $\rho$. 

We endow $M_{\rho}(X)$ with the coarsest
$\sigma$-algebra such that, for all Borel sets $B \subset
X$, the mapping $M_\rho(X) \ni \nu \mapsto
\nu(B)$ is Borel.
Let $\mathcal{S}$ denote the $\sigma$-algebra of $X$, and let
 $\mathcal{S}^H$ be the subalgebra of $H$-invariant Borel sets.
Clearly $\mathcal{S}^H$ is a subalgebra of $\mathcal{S}$. 
% Given $\nu \in M_{\rho}(X)$ and a positive Borel function $f$,
% denote its image measure under the canonical map by $\overline{\nu}$.
The conditional expectation of $f$ with respect to $\nu \in M_\rho(X)$ is a Borel function
\[
 E_{\nu}(f|\mathcal{S}^H) : X \to \mathbb{R}_0^+
\] which is $H$-invariant and fulfills 
\[
 \int_B f(\xi) d\nu(\xi) = \int_{B}  E_{\nu}(f|\mathcal{S}^H)(\xi) d\nu(\xi)~,
\] for all $H$-invariant Borel sets $B$. The conditional expectations always exists
and is $\nu$-a.e. unique \cite[5.1.15]{Str} 

By \cite[Theorem 5.2]{GrSch}, there exists an $H$-invariant map $p : X \to M_{\rho}(X),
x \mapsto p_\xi$ such that $p_\xi$ is ergodic, and in addition, for every $\nu \in M_{\rho}(X)$ 
and every positive Borel function $f$ on $X$,
\begin{equation} \label{eqn:ergdc_cond_exp}
 E_{\nu}(f|\mathcal{S}^H) (\xi) = \int_X f(\omega) dp_\xi(\omega)~.
\end{equation} holds for $\nu$-almost all $\xi \in X$.

\begin{lemma}
Assume that $(\overline{\lambda},(\beta_{\mathcal{O}})_{\mathcal{O}
\subset X})$ is a decomposition of $\lambda$ into probability measures over the orbits. Let $\rho$
be the cocycle of $\lambda$, and let $p: X \to M_{\rho}(X)$
denote the ergodic decomposition associated to $\rho$.
 There exists a conull, $H$-invariant Borel set $Y \subset X$ such that  $p_\xi = \beta_{\xi.H}$, for all $\xi \in Y$.
\end{lemma}

\begin{prf}
 We first observe that for almost all $\mathcal{O}$,
$\beta_{\mathcal{O}} \in M_{\rho}(X)$: For Borel subsets $B \subset X$ and $H$-invariant $C \subset X$,
\begin{eqnarray*}
 \int_{C} \beta_{\mathcal{O}} (B.h) ~d\overline{\lambda}(\mathcal{O})& = &  \lambda(B.h \cap C)  \\
 & = & \int_{C} \mathbf{1}_B (\xi) e^{\rho(\xi,h)} ~d\lambda(\xi) \\
 & = & \int_{C} \int_B e^{\rho(\xi,h)} ~d\beta_{\mathcal{O}}(\xi) ~d\overline{\lambda}(\mathcal{O}) ~,
\end{eqnarray*} and thus, for all $h \in H$,
\begin{equation} \label{eqn:cccle_beta}
 \beta_{\mathcal{O}} (B.h) = \int_B e^{\rho(\xi,h)} ~d\beta_{\mathcal{O}}(\xi) ~,
\end{equation} valid for a  $\overline{\lambda}$-conull set of orbits $\mathcal{O}$ that may still depend on $h \in H$ and $B$.

By Fubini's theorem, for each $B \in \mathcal{S}$ there exists $Y(B) \subset X$ Borel, $H$-invariant and conull such that 
 (\ref{eqn:cccle_beta}) holds for all orbits $\mathcal{O} \subset Y(B)$ and all $h \in T(\mathcal{O},B)$, with 
$T(\mathcal{O},B) \subset H$ Borel, conull.
Next pick a generating sequence $(B_k)_{k \in \mathbb{N}}$ of $\mathcal{S}$, and define
\[
 Y = \bigcap_{k \in \mathbb{N}} Y(B_k)~,~\forall \mathcal{O} \subset Y~:~T(\mathcal{O}) = 
\bigcap_{k \in \mathbb{N}} T(\mathcal{O},B_k)~.
\] Then (\ref{eqn:cccle_beta}) holds for all $B \in \mathcal{S}$, $ \mathcal{O} \subset Y$ and
$h \in T(\mathcal{O})$, since both sides of (\ref{eqn:cccle_beta}) define a Borel measure, hence
coincide on a $\sigma$-algebra. 

Now fix $\mathcal{O} \subset Y$, and define 
\[
 H(\mathcal{O}) = \{ h \in H : \forall B \in \mathcal{S} : \beta_{\mathcal{O}} (B.h) = \int_B e^{\rho(\xi,h)} ~d\beta_{\mathcal{O}}(\xi) \}~.
\]
We claim that $H(\mathcal{O})$ is a subgroup of $H$: Assume that $h \in H(\mathcal{O})$.
Then  (\ref{eqn:cccle_beta})  extends to positive Borel functions $f$, yielding
\begin{equation} \label{eqn:ccle_beta_fun}
\int_X f(\xi.h^{-1}) d\beta_{\mathcal{O}}(\xi) = \int_X f(\xi) e^{\rho(\xi,h)} d\beta_{\mathcal{O}}(\xi) ~.
\end{equation}
Furthermore, the cocycle properties  (\ref{eqn:cc_cond}) and  (\ref{eqn:cc_fix}) entail that
$\rho(\xi,h^{-1}) = -\rho(\xi.h^{-1},h)$. 
Using this, we can compute
\begin{eqnarray*}
 \int_B e^{\rho(\xi,h^{-1})} d\beta_{\mathcal{O}}(\xi)  & = & 
 \int_B e^{-\rho(\xi.h^{-1},h)} d\beta_{\mathcal{O}}(\xi) \\
 & = & \int_X \mathbf{1}_B(\xi)  e^{\rho(\xi.h^{-1},h)} d\beta_{\mathcal{O}}(\xi) \\
  & \stackrel{(\ref{eqn:ccle_beta_fun})}{=} & \int_X \mathbf{1}_B(\xi.h) e^{-\rho(\xi,h)} e^{\rho(\xi,h)} d\beta_{\mathcal{O}}(\xi) \\
 & = & \int_{B.h^{-1}} d\beta_{\mathcal{O}}(\xi) \\
 & = & \beta_{\mathcal{O}}(B.h^{-1}) 
\end{eqnarray*}
which proves $h^{-1} \in H(\mathcal{O})$. 

Next let $g,h \in H(\mathcal{O})$. Then, since $g \in \mathcal{O}$,
\begin{eqnarray*}
 \beta_{\mathcal{O}}(B.hg) & = & \int_{B.h} e^{\rho(\xi,g)} d\beta_{\mathcal{O}}(\xi) \\
  & = & \int_X \mathbf{1}_{B}(\xi.h^{-1}) e^{\rho(\xi,g)} d\beta_{\mathcal{O}}(\xi) \\
 &\stackrel{(\ref{eqn:ccle_beta_fun})}{=} & \int_X \mathbf{1}_B (\xi) e^{\rho(\xi,h)}
e^{\rho(\xi.h,g)}  d\beta_{\mathcal{O}}(\xi) \\ 
 &\stackrel{(\ref{eqn:cc_cond})}{=} & \int_X \mathbf{1}_B (\xi) e^{\rho(\xi,hg)} d\beta_{\mathcal{O}}(\xi) ~,
\end{eqnarray*}
and therefore $hg \in  H(\mathcal{O})$.

Hence $H(\mathcal{O}) \subset H$ is a subgroup, with $H(\mathcal{O}) \supset T(\mathcal{O})$. 
In particular $H(\mathcal{O}) \supset T(\mathcal{O}) T(\mathcal{O})^{-1}$, and since $T(\mathcal{O})$
has positive Haar measure, $H(\mathcal{O})$ contains a nonempty open subset \cite[Proposition III.12.3]{FeDo}. 
Hence $H(\mathcal{O})$ is an open subgroup, and therefore closed. On the other hand, 
$T(\mathcal{O})$ is conull and thus dense in $H$, whence finally $H = H(\mathcal{O})$. But this shows 
$\beta_{\mathcal{O}} \in M_\rho(X)$ for all $\mathcal{O} \subset Y$. 

Then, since $\beta_{\mathcal{O}}$ is supported in $\mathcal{O}$, it follows for every nonnegative Borel
functions $f$ and $\xi \in \mathcal{O}$
that 
 \[
  \int_X f(\omega) d\beta_{\mathcal{O}} (\omega) =  E_{\beta_{\mathcal{O}}}(f|\mathcal{S}^H) (\xi) =
 \int_X f(\omega) dp_\xi(\omega)~,
\] 
where the second equation is due to (\ref{eqn:ergdc_cond_exp}). But this means that $ \beta_\mathcal{O} = p_\xi$. 
\end{prf}

Hence, after passing to a suitable conull $H$-invariant subset, we may assume that $\beta_{\xi.H} = p_\xi$
holds for all $\xi \in X$. In particular, we may assume in the following that $p$ separates the orbits in $X$.

Denote by  $\mathcal{T}$ the coarsest $\sigma$-algebra on $X$ making $p$ a Borel map.
Since the $\beta_{\mathcal{O}}$ are a measurable family, $p: X \to M_{\rho}(X)$ is
clearly Borel, thus $\mathcal{T} \subset \mathcal{S}$.
On the other hand, by \cite[Theorem 5.2]{GrSch} $\mathcal{T}$ is
countably generated.

Since $p$ is $H$-invariant, the elements of $\mathcal{T}$ are
$H$-invariant as well. Hence $q: X \to
X/H$ induces an isomorphism of $\sigma$-algebras
between $\mathcal{T}$ and its image $\overline{\mathcal{T}} = \{
q(A): A \in \mathcal{T} \}$. In particular, the latter is countably
generated as well, and it is contained in the quotient
$\sigma$-algebra on $X/H$. Furthermore, it is clearly separated, since $p$
separates the orbits. But then the quotient $\sigma$-algebra, being finer than $\overline{\mathcal{T}}$,
is countably separated. Hence, by \cite[Proposition 2.9]{AuMo}, it follows that $X/H$
is an analytic Borel space. But then there exists a conull Borel subset $A \subset X/H$
which is standard (see \cite{AuMo}, remarks following I.2.13 ). This shows (c). 

\subsection{Proof of $(c) \Rightarrow (d) \Rightarrow (b)$}

For $(c) \Rightarrow (d)$ we may assume, after passing to a suitable
conull subset, that $X/H$ is standard. Then \cite[Proposition 2.15]{AuMo} yields
a $\overline{\lambda}$-conull set $V \subset X/H$
and a Borel cross-section $\sigma: V \to q^{-1}(V)$. 
Then $\sigma$ is injective, and $V$ is standard, as a Borel subset of $X/H$. 
But then $\sigma(V)$ is Borel, by \cite[Proposition 2.5]{AuMo}, and it meets every 
orbit contained in $V$ in precisely one point.  

Finally, $(d) \Rightarrow (b)$ follows by \cite[Lemma 11.1]{Ma}.

\subsection{Proof of $(d) \Rightarrow (a)$}

Now assume $(d)$, and that all stabilizers are compact.  
Let $Y \subset X$ be $H$-invariant and conull, and let $C \subset Y$
be a Borel transversal for the orbits in $Y$. Let $K \subset H$
denote a compact neighborhood of the identity, and $V = C.K = \{
\xi.h: \xi \in  C, h \in K \}$. Then $V$ is an analytic subset of $X$, as
the Borel image of the standard set $C \times K$ in the countably
generated space $Y$. Since analytic sets are universally measurable 
(confer \cite{AuMo}, page 11), $V$ is $\lambda$-measurable.
Hence there exist sets $U \subset V \subset W$, with $U$, $W$ Borel
and $\lambda(W \setminus U) = 0$.

We intend to use $\varphi = \mathbf{1}_W$ to show weak admissibility. 
This amounts to showing, for almost all $\xi \in  X$, that 
\begin{equation} \label{eqn:adm_chiW} 0 < \mu_{\mathcal{O}}(W) = \mu_H(\{ h : \xi.h \in W \}) < \infty ~,\end{equation}
for $\mathcal{O} \ni \xi$ .
In order
to do this, we first consider $\mathbf{1}_V$. Note that for every $\xi \in  X$
\[
 \{ h \in H : \xi.h \in V \} = H_\xi K
\] is compact. Since the canonical map $H_\xi \setminus H \to \mathcal{O}$ is a Borel
isomorphism, it follows that $V \cap \xi.H$ is in fact a Borel set. In addition, since
$H_\xi K$ is a compact neighborhood of the identity element, 
\begin{equation} \label{eqn:adm_chiV} 0 <
\mu_{\mathcal{O}}(V \cap \xi.H) < \infty ~.
\end{equation}

In order to conclude (\ref{eqn:adm_chiW}) from this, we use $(d) \Rightarrow (b)$ and decompose $\lambda$
into a family $(\beta_{\mathcal{O}})_{\mathcal{O} \subset X}$ of measures supported on the orbits. Then almost every $\beta_{\mathcal{O}}$
is equivalent to a finite quasi-invariant measure $\tilde{\beta}_{\mathcal{O}}$. With respect to the topology induced
by the canonical bijection $H_\xi \setminus H \to \mathcal{O}$, the finite measure $\tilde{\beta}_{\mathcal{O}}$ becomes
regular \cite[Theorem 7.8]{Fo_RA}. On the other hand, $\mu_{\mathcal{O}}$ is also a regular quasi-invariant measure,
hence $\mu_{\mathcal{O}}$ is equivalent to $\tilde{\beta}_{\mathcal{O}}$ by \cite[14.9]{FeDo}, and thus finally
to $\beta_{\mathcal{O}}$.
% 
% By definition
% \[
% U \cap \xi.H \subset V \cap \xi.H \subset W \cap \xi.H~.
% \] 

Now $\lambda(W \setminus U) = 0$ entails
$\beta_{\mathcal{O}}(W \setminus U) = 0$, for almost all orbits $\mathcal{O}$.
Since $\beta_{\mathcal{O}}$ is equivalent to $\mu_{\mathcal{O}}$, it follows
for these orbits that
\[ \mu_{\mathcal{O}}((W \cap \xi.H) \setminus (V \cap \xi.H)) \le
\mu_{\mathcal{O}}((W \cap \xi.H) \setminus (U \cap \xi.H)) = 0 ~,\] 
with $\xi \in \mathcal{O}$. Thus $\mu_{\mathcal{O}}(W) = \mu_{\mathcal{O}}(W \cap \xi.H) = \mu_{\mathcal{O}} (V \cap \xi.H)$, and thus  (\ref{eqn:adm_chiV}) implies  
(\ref{eqn:adm_chiW}).

\section{Admissible vectors versus weakly admissible vectors}

Throughout this section, $X \subset \widehat{A}$ is Borel-measurable, $H$-invariant,
and $\mathcal{H} = \mathcal{H}_X$. For explicit reference to the results of the previous
two sections, let $\lambda$ denote Haar measure on $\widehat{A}$. 

We assume the existence of a weakly admissible
vector in $\mathcal{H}$, and want to clarify which additional criteria must be met to ensure the existence
of an admissible vector. 

The main tool for this purpose will be the decomposition of Haar measure on $X$.
The discussion in this section closely follows \cite[Section 5.2]{Fu_LN}, but we have chosen
to spell out most details for two reasons: First, we start  
from somewhat more general assumptions, and secondly, the arguments in \cite[Section 5.2]{Fu_LN},
are partly flawed. This applies in particular to \cite[Lemma 5.9]{Fu_LN}, which is an analog of the 
following result. Thus the following serves both as erratum and generalization to some of the
results in \cite{Fu_LN}.

\begin{lemma} \label{lem:meas_decomp_Planch}
Assume that $\mathcal{H}$ has a weakly admissible vector. 
\begin{enumerate}
 \item[(a)] Fix any pseudo-image $\overline{\lambda}$ of Plancherel measure on $X$. There exists an essentially unique
family of measures $(\beta_{\mathcal{O}})_{\mathcal{O} \subset X}$
such that $d\xi = d\beta_{\mathcal{O}}(\xi) d\overline{\lambda}(\mathcal{O})$.
 \item[(b)] For every orbit $\mathcal{O} \subset X$ let $\mu_{\mathcal{O}}$ be as in
Lemma \ref{lem:defnmu}. There exists an essentially unique Borel function $\kappa: X \to \mathbb{R}_0^+$
such that, for $\overline{\lambda}$-almost all orbits
\[
 \frac{d\beta_{\mathcal{O}}}{d\mu_{\mathcal{O}}}(\xi) = \kappa (\xi)~.
\]
\item[(c)] $\kappa$ can be chosen in such a way that for all $h \in H$ and all $\xi$
in a fixed $H$-invariant conull set, 
\[
 \kappa(\xi.h) =  \kappa(\xi) \Delta_G(0,h)^{-1} .
\]
In particular, $\kappa$ is $H$-invariant iff $G$ is unimodular. In this case, $\lambda$ has a
decomposition $(\overline{\lambda},
(\mu_{\mathcal{O}})_{\mathcal{O} \subset X})$  of $\lambda$, where $\overline{\lambda}$ 
is a suitable $\sigma$-finite measure. 
\end{enumerate}
\end{lemma}
\begin{prf}
Part (a) is Theorem \ref{thm:main} $(a) \Rightarrow (b)$. For part (b) let $\mu$ be as defined in Lemma
\ref{lem:defnmu}. Then $\mu$ and $\lambda$ are equivalent $\sigma$-finite measures, as was shown in 
the proof of \ref{thm:main} $(a) \Rightarrow (b)$, and we find that 
\[ \kappa(\xi) = \frac{d\lambda}{d\mu}(\xi) = \frac{d\beta_{\mathcal{O}}}{d\mu_{\mathcal{O}}}(\xi) \] 
is the desired global Radon-Nikodym-derivative. Thus
(b) follows. For part (c), we let $\mu_h(B) = \mu(B.h)$, and $\lambda_h(B) = \lambda(B.h)$.
Then 
\[
\frac{d\mu_h}{d\mu} (\xi) = \Delta_H(h) ~,~\frac{d\lambda_h}{d\lambda}(\xi) = \delta(h)~.
\] For any nonnegative Borel map $f$ on $X$, the definition of $\mu_h$ entails 
\[
 \int_X f(\xi) d\mu_h(\xi) = \int_X f(\xi.h^{-1}) d\mu(\xi)~.
\]
It follows for $h \in H$ and arbitrary Borel sets $B \subset X$, that 
\begin{eqnarray*}
 \int_B \frac{d \lambda}{d\mu}(\xi.h) ~d\mu_h(\xi) & = & \int_X \mathbf{1}_B (\xi) \frac{d \lambda}{d\mu}(\xi.h) ~d\mu_h(\xi) \\
 & = & \int_X \mathbf{1}_B (\xi.h^{-1}) \frac{d \lambda}{d\mu}(\xi) ~d\mu(\xi) \\
 & = & \int_{B.h}  \frac{d \lambda}{d\mu}(\xi) ~d\mu(\xi)  \\ & = &  \lambda_h(B) \\ 
 & = & \int_B \frac{d\lambda_h}{d\mu_h}(\xi)~ d\mu_h(\xi)~,
\end{eqnarray*}
and thus 
\[
  \frac{d\lambda_h}{d\mu_h}(\xi) = \frac{d\lambda}{d\mu}(\xi.h) \hspace{1.5cm} (\lambda-\mbox{a.e. })
\]
Hence, for a.e. $\xi \in X$, the chain rule for Radon-Nikodym-derivatives yields
\[
 \kappa(\xi.h) = \frac{d\lambda}{d\mu}(\xi.h) = \frac{d\lambda_h}{d\mu_h}(\xi) =
\frac{d\lambda_h}{d\lambda}(\xi) \frac{d\lambda}{d\mu}(\xi) \frac{d\mu}{d\mu_h}(\xi) = \kappa(\xi) \frac{\delta(h)}{\Delta_H(h)}~,
\] which is the desired equality, except that the conull subset of $X$ on which it holds may
still depend on $h$. However, by \cite[B.5]{Zi}, one finds a conull invariant Borel subset
of $X$ on which the relation holds everywhere, independent of $h$. 

If  $\kappa$ is constant on the orbits,
it defines a Borel mapping $\overline{\kappa}$ on
$X/H$. Replacing each $\beta_{\mathcal{O}}$ by
$\mu_{\mathcal{O}}$, we can make up for it by taking
$\overline{\kappa}(\mathcal{O}) d\overline{\lambda}(\mathcal{O})$ as the new
measure on the orbit space. The result is a $\sigma$-finite measure $\overline{\kappa} d\overline{\lambda}$.
\end{prf}

The next result clarifies the role of the specific choice of $\overline{\lambda}$. 
\begin{thm} \label{thm:main_ex_adm}
Let $\mathcal{H} = \mathcal{H}_X \subset {\rm L}^2(A)$ be closed and $\pi$-invariant. 

There exists an admissible vector for $\mathcal{H}$ iff  there
exists a weakly admissible vector, and in addition, 
\begin{enumerate}
 \item $G$ is nonunimodular; or 
 \item $G$ is unimodular, and with $\overline{\lambda}$ chosen according to Lemma \ref{lem:meas_decomp_Planch} (c): 
 \[ \overline{\lambda}(X/H) < \infty ~.\] 
\end{enumerate}
\end{thm}
\begin{prf}
First assume that $G$ is unimodular, and that $\psi$ is an admissible vector. 
Then the Plancherel theorem and the measure decomposition 
over the orbits, with $\overline{\lambda}$ as in Lemma \ref{lem:meas_decomp_Planch} (c) and $\beta_{\mathcal{O}}
= \mu_{\mathcal{O}}$, allows to compute
\begin{eqnarray*}
\| \psi \|_2^2 & = & \int_{X} |\widehat{\psi}(\xi)|^2~d\lambda(\xi)  \\
& = & \int_{X/H} \int_{H} |\widehat{\psi}(\xi)|^2~d\mu_{\mathcal{O}}(\xi) ~d\overline{\lambda}(\mathcal{O}) \\
& = & \int_{X/H} \int_{H} |\widehat{\psi}(\xi.h)|^2~dh ~d\overline{\lambda}(\mathcal{O}) \\
& = & \int_{X/H} 1 ~d\overline{\lambda}(\mathcal{O}) \\
& = & \overline{\lambda}(X/H)~.
\end{eqnarray*} Here the penultimate equality was due to admissibility of $\psi$. 
In particular $\overline{\lambda}(X/H) < \infty$. 

For the converse, assume that $\psi_0$ is a weakly admissible vector. Define 
\[
 \Phi(\xi) = (\int_{H} |\widehat{\psi_0}(\xi.h)|^2~dh)^{1/2}~.
\] By assumption, $0< \Phi(\xi) < 1$ a.e.
Let $\varphi(\xi) = \widehat{\psi_0}(\xi)/\Phi(\xi)$.
It follows that 
\[
 \int_{H} |\varphi(\xi.h)|^2 d\xi = 1~.
\] 

If $G$ is unimodular, the measure decomposition allows to compute
\[
 \| \varphi \|_2^2 = \int_{X/H} \int_{H} |\varphi(\xi.h)|^2 ~d\xi d\overline{\lambda}(\mathcal{O}) = 
\overline{\lambda}(X/H) ~.
\] Thus, if $\overline{\lambda}(X/H)< \infty$, the inverse Plancherel transform of $\varphi$ is
admissible for $\mathcal{H}_X$. 

Finally,  assume that $G$ is nonunimodular.  Then $\Delta_G$ is nontrivial on $H$, and there
exists $h_0 \in H$ such that $\Delta_G(h_0) < 1/2$. Since $\overline{\lambda}$ is $\sigma$-finite, 
we can write $X$ as a disjoint union $X = \bigcup_{n\in \mathbb{N}} V_n$, where $V_n \subset
X$ Borel, $H$-invariant and with $\overline{\lambda}(V_n/H)< \infty$. 
Since $\widehat{\psi} \in {\rm L}^2(\widehat{A})$, 
\[
 \Psi : \xi \mapsto \left( \int_{\xi.H} |\varphi(\xi)|^2 d\beta_{\xi.H}(\xi) \right)^{1/2}
\] is finite a.e., and  we may in addition
assume that $\Psi$ is bounded on each $V_n$; in particular,
the functions $( \mathbf{1}_{V_n} \cdot \Psi)_{n \in \mathbb{N}}$ 
are square-integrable.

Now pick a sequence $(k_n)_{n \in \mathbb{N}}$ of integers satisfying 
\[
 2^{-k_n} \|  \mathbf{1}_{V_n} \cdot \Psi  \|_2^2 < 2^{-n}~,
\] and let 
\[
 \nu(\xi) = \sum_{n \in \mathbb{N}} \Delta_H(h_0)^{k_n} \varphi(\xi.h_0^{k_n})~.
\]
On the one hand, 
\begin{eqnarray*}
 \int_{X} |\nu(\xi)|^2 d\lambda(\xi) & = & \int_{X/H} \int_{\mathcal{O}} |\nu(\xi)|^2 d\beta_{\mathcal{O}}(\xi) 
d \overline{\lambda}(\mathcal{O}) \\
 & = &  \sum_{n \in \mathbb{N}} \int_{V_n} \int_{\mathcal{O}} \Delta_H(h_0)^{k_n} |\varphi(\xi.h_0^{k_n})|^2
d\beta_{\mathcal{O}}(\xi) d \overline{\lambda}(\mathcal{O}) \\
 & = & \sum_{n \in \mathbb{N}} \int_{V_n} \int_{\mathcal{O}} \Delta_H(h_0)^{k_n} \delta(h_0)^{-k_n} |\varphi(\xi)|^2
d\beta_{\mathcal{O}}(\xi) d \overline{\lambda}(\mathcal{O}) \\
& = & \sum_{n \in \mathbb{N}} \Delta_G(h_0)^{k_n} \int_{V_n} \int_{\mathcal{O}} |\varphi(\xi)|^2
d\beta_{\mathcal{O}}(\xi) d \overline{\lambda}(\mathcal{O}) \\
& = & \sum_{n \in \mathbb{N}}  \Delta_G(h_0)^{k_n}  \|  \mathbf{1}_{V_n} \cdot \Psi  \|_2^2 \\
& \le & \sum_{n \in \mathbb{N}}  2^{-k_n}  \|  \mathbf{1}_{V_n} \cdot \Psi  \|_2^2 \\
& < &  \infty,
\end{eqnarray*} by choice of the $k_n$.
Hence $\nu$ is square-integrable. Moreover, the Calder\'on condition
is also easily verified: For $x\in V_n$,
\begin{eqnarray*}
 \int_H |\nu(\xi.h)|^2 dh & = & \int_H |\varphi(\xi.h h_0^{k_n})|^2 \Delta_H(h_0)^{k_n} dh \\
 & = & \int_H |\varphi(\xi.h )|^2 dh \\
 & = & 1~,
\end{eqnarray*}
by construction of $\varphi$.
Thus the inverse Plancherel transform of $\nu$ is the desired admissible vector. 
\end{prf}

% 
% 
% \begin{cor}
%  There exists an admissible vector for ${\rm L}^2(A)$ iff there exists a weakly admissible vector for ${\rm L}^2(A)$
% and in addition, $G$ is nonunimodular or $G$ is discrete. 
% \end{cor}
% \begin{prf}
%  The nonunimodular part is immediate from the theorem. The unimodular part follows from the
% additional observation that $\widehat{A}$ has finite measure iff $\widehat{A}$ is compact, 
% which is in turn equivalent to discreteness of $A$. 
% \end{prf}

\begin{rem}
For unimodular semidirect products, we do not have a clean-cut and complete characterization of the group having an 
admissible vector for all of ${\rm L}^2(A)$. A straightforward adaptation of the proof for \cite[Proposition 5.14]{Fu_LN}
allows to describe a rather general setting in which ${\rm L}^2(A)$ does not have an admissible vector: 

Suppose that $G = A \rtimes H$ is unimodular, and has a weakly admissible vectors. 
Let $r$ be a topological automorphism of $A$. 
We assume that $r$ has the following properties:
\begin{enumerate}
 \item[(i)] $r$ normalizes $H$.
 \item[(ii)] For any (hence all) $B \subset H$  and $C \subset A$ of positive finite Haar measure, 
\[
\frac{|rBr^{-1}|}{|B|} \not= \frac{|r(C)|}{|C|}~.
\] 
\end{enumerate}
Then $\overline{\lambda}(X/H) = \infty$.  In particular, there exists no admissible
vector for ${\rm L}^2(A)$. 

This result applies in particular to $A = \mathbb{R}^d$:  Choose
$r = s \cdot {\rm Id}_{\mathbb{R}^d}$, with $s \not=1$. Then $r$ commutes
with all elements of the matrix group $H$. In particular,
conjugation with $r$ leaves Haar measure on $H$ invariant, whereas $\frac{|r(C)|}{|C|} = s^d$.
Thus (ii) is ensured, which proves that there exist no admissible vectors in this case. 
\end{rem}

The final result concerns irreducible representations. Recall that irreducible
representations with admissible vectors are called discrete series representations. 
Most early sources generalizing wavelets to higher dimension restricted their attention
to the discrete series case, e.g. \cite{Mu,BeTa,Fu}. The implication
$(b) \Rightarrow (a)$ of the following result has been proved for $A=\mathbb{R}^n$ in \cite{Fu}.
However, the converse was previously known only for 
$A = \mathbb{R}^n$ and $H \subset {\rm GL}(n,\mathbb{R})$ discrete, where it boils
down to stating that no discrete series representation of that type exists, see \cite[Remark 12]{Fu}.

\begin{cor} \label{cor:discrete_series}
 Let $\mathcal{H}_X \subset {\rm L}^2(A)$ be a nontrivial closed $\pi$-invariant subspace. The following
are equivalent:
\begin{enumerate}
 \item[(a)] The restriction of $\pi$ to $\mathcal{H}_X$ is a discrete series representation.
 \item[(b)]  There exists an orbit $\mathcal{O} \subset X$ such that $|X \setminus \mathcal{O}| =0$,
with associated compact stabilizers. 
\end{enumerate}
\end{cor}
\begin{prf}
For $(b) \Rightarrow (a)$, the arguments given in \cite{Fu} immediately carry over;
see also \cite{Ka_etal}.

Conversely, assume that $\pi$ restricted to $\mathcal{H}_X$ is in the discrete series.  
If $X = W \cup V$ with disjoint, $H$-invariant Borel sets $U,W$ of positive measure, then
$\mathcal{H}_X = \mathcal{H}_W \oplus \mathcal{H}_V$ contradicts irreducibility.
Thus the action of $H$ on $X$ is ergodic with respect to Haar measure. Since $\overline{\lambda}$ is standard
on $X/H$, it follows by \cite[Chapter I, Proposition 3.9]{AuMo} that there
exists a conull orbit. The associated stabilizers must be compact by Theorem \ref{thm:main_grp}.
\end{prf}

\begin{rem}
 The measure decompositions discussed in this paper are closely related to direct integral theory. 
In order to see this connection, first note that the quasi-regular representation $\pi$ is type I: Its
commuting algebra is contained in the commuting algebra of the regular representation of $A$
on ${\rm L}^2(A)$; $A$ being abelian, the latter algebra is commutative. Hence $\pi$
is multiplicity-free, in particular type I. It therefore has a unique direct integral decomposition into irreducibles, which
is closely related to the ergodic decomposition of $\lambda$.

For the sake of simplicity, let us assume that there exists a weakly admissible vector, so that the ergodic
decomposition is in fact a decomposition over the orbits. Then the measure decomposition $d\lambda(\xi) =
 d\beta_{\mathcal{O}}(\xi)  d\overline{\lambda} (\mathcal{O})$ gives rise to a direct integral decomposition
\[
 {\rm L}^2(\widehat{A}) \simeq \int_{X/H}^\oplus {\rm L}^2(\mathcal{O},d\beta_{\mathcal{O}}) d\overline{\lambda}(A)~.
\]
It can be shown that this decomposition also applies to the representation, yielding
\[
 \pi \simeq \int_{X/H}^\oplus {\rm Ind}_{A \semdir H_\xi}^G (\xi \times \mathbf{1})~ d\overline{\lambda}(\mathcal{O})~, 
\]
where $\mathbf{1}$ denotes the trivial representation of $H_\xi$. 
By Mackey's theory, the induced representations are irreducible (and pairwise inequivalent), thus we have decomposed $\pi$
into irreducibles.

But the orbit space $\widehat{A}/H$ also occurs in the direct integral decomposition of the 
regular representation of $G$. In fact, 
the existence of a weakly admissible vector for ${\rm L}^2(A)$ implies that the regular representation of $G$
is type I: By Theorem \ref{thm:main_grp}, almost all stabilizers are compact, and the dual orbit space is standard
up to a set of measure zero. Note that compactness of the stabilizer $H_\xi$ entails 
that $H_{\xi}$ has a type I regular $\omega$-representation, where $\omega$ denotes an arbitrary multiplier on $H_\xi$.
Furthermore, the orbit space is standard (outside a set of measure zero). Thus, by \cite[Theorem 2.3]{KlLi}, it follows
that the regular representation of $G$ is type I, and that the Plancherel measure of $G$ is obtained as
fibred measure with base space given by $\widehat{A}/H$, base measure given by $\overline{\lambda}$, 
and fibres given by the $\omega_\xi$-duals of the $H_\xi$, where $\omega_\xi$ are suitably chosen
multipliers on $H_\xi$. 

Now the connection between $\pi$ and the left regular representation can also be realized by observing
that Mackey's construction yields a mapping
\[
 \widehat{A}/H \ni \xi.H \mapsto {\rm Ind}_{A \semdir H_\xi}^G (\xi \times \mathbf{1}) \in \widehat{G}
\] identifying $\widehat{A}/H$ with a (Borel) subset of $\widehat{G}$. 
It then becomes apparent that the
measure $\overline{\lambda}$ underlying the direct integral decomposition of $\pi$ is nothing
but the restriction of Plancherel measure of $G$ to this subset. This is an alternative proof for the
containment of $\pi$ in the regular representation. 
This type of reasoning, using direct integral decompositions to study existence of inversion
formulae, has been developed systematically in \cite{Fu_LN}.  In particular, \cite[Section 5.3]{Fu_LN} 
contains a rigourous investigation of the double role of the measure $\overline{\lambda}$ 
in decomposing both $\pi$ and the regular representation. 
Note however that the underlying assumption of \cite{Fu_LN} is that $G$ is type I.
By contrast, we make no such initial assumption on $G$, and obtain that  
the regular representation is type I as a consequence of the existence of weakly admissible
vectors. 
\end{rem}

\begin{rem}
 The results presented in this paper are satisfactory to a certain degree, since they provide
a sharp characterization. However, we are not aware of an easy general procedure for the explicit verification 
of the criteria in concrete cases. Also, we do not know how our characterization relates to 
other criteria, in particular compactness of almost all $\epsilon$-stabilizers, proven to be sufficient
in \cite{LWWW}.

To our knowledge, the first systematic and substantial investigation of regularity properties for
orbit spaces was carried out by Glimm \cite{Gl}, who proved that  standardness of the orbit space
of a second countable locally compact group $H$ acting continuously on a second countable locally compact space $X$
is equivalent to a variety of conditions, most notably countable separatedness of $X/H$,
or the existence of a Borel cross-section, or local compactness of the orbits in the relative topology. 
On the one hand, these results closely resemble our conditions
$(c)$ and $(d)$ from above, but also the $\epsilon$-stabilizer condition: To see this, note that 
compactness of $H_{\epsilon,\xi}$, for some $\epsilon>0$,
is equivalent to $i)$ compactness of $H_\xi$, and in addition $ii)$, local compactness of the 
orbit $\xi.H$ in the relative topology (cf. the proof of \cite[Proposition 5.7]{Fu_LN}). 
Hence, if Glimm's results were applicable to our setting, they would imply
that weak admissibility of the dual action is equivalent to existence of a compact $\epsilon$-stabilizer,
for a.e. $\xi$. 

However, a direct application of Glimm's results to our setting is 
impeded by the fact that by definition, weak admissibility only concerns the behaviour of the orbits
in a suitable conull subset. In particular, weak admissibility is robust under passage to a conull
invariant subset, whereas the assumptions underlying Glimm's characterization can be seriously
affected by this step:  A conull Borel subset of a locally compact space no longer needs to 
be locally compact. It was mostly this obstacle that stopped previous efforts of the author to characterize 
weakly admissible group actions. Attempts to use more recent generalizations 
of Glimm's results for the study of admissibility got stuck for similar reasons.
\end{rem}

\begin{rem}
Throughout this paper, all groups have been assumed to be second countable. Most of the 
measure-theoretic arguments in this paper strongly rely on countability assumptions, and it is currently
open to what extent our results can be generalized beyond second countable groups. 
\end{rem}

\section*{Acknowledgements}
Part of the research for this paper was carried out at the Universitat Aut\'onoma de Barcelona,
Departament de Matem\'atiques, and I am grateful to both institutions for their hospitality.
I would specifically like to thank Joaquim Bruna for interesting and fruitful discussions.
Thanks are also due to G.~Greschonig for clarifying comments on \cite{GrSch}.

\end{document}